\documentclass[a4j,12pt, leqno]{article}

\usepackage{amsmath,amsthm,amssymb}
\usepackage[dvipdfmx]{color}
\usepackage[dvipdfmx]{graphicx}
\usepackage{latexsym}
\usepackage{ascmac}
\usepackage{longtable}
\usepackage{comment}
\usepackage{color}
\usepackage{setspace}

\newcommand{\argmax}{\mathop{\rm argmax}\limits}
\newcommand{\sign}{\operatorname{sign}} 
\newcommand{\mb}{\mathbb}

\theoremstyle{plain}
\newtheorem{thm}{Theorem}
\newtheorem{lem}[thm]{Lemma}
\newtheorem{cor}{Corollary}
\newtheorem{prop}[thm]{Proposition}
\newtheorem{asm}{Assumption}

\newtheorem{exm}{Example}
\newtheorem{defn}{Definition}

\theoremstyle{remark}

\newtheorem{prf}{Proof}

\begin{document}

\title{Duality of nonconvex optimization with\\ positively homogeneous functions%
\thanks{
This work 
was supported in part by a Grant-in-Aid for Scientific Research (C) (17K00032) from 
Japan Society for the Promotion of Science.
}}

\author{
Shota Yamanaka%
\thanks{Department of Applied Mathematics and Physics, 
Graduate School of Informatics, Kyoto University,
Kyoto 606-8501, Japan (\texttt{shota@amp.i.kyoto-u.ac.jp, nobuo@i.kyoto-u.ac.jp}).}
\ and
Nobuo Yamashita$^\dagger$
}

\maketitle

\begin{abstract}
\noindent
We consider an optimization problem with positively homogeneous functions
in its objective and constraint functions.
Examples of such positively homogeneous functions include
the absolute value function and the $p$-norm function,
where $p$ is a positive real number.
The problem, which is not necessarily convex, extends the absolute value optimization
proposed in [O.~L.~Mangasarian, Absolute value programming,
Computational Optimization and Applications 36 (2007) pp. 43--53].
In this work, we propose a dual formulation that, differently from the Lagrangian dual approach,
has a closed-form and some interesting properties.
In particular, we discuss the relation between the Lagrangian duality and the one proposed here,
and give some sufficient conditions under which these dual problems coincide.
Finally, we show that some well-known problems, e.g.,
sum of norms optimization and the group Lasso-type optimization problems,
can be reformulated as positively homogeneous optimization problems.\\

\end{abstract}

\section{Introduction}

Recently, the so-called \textit{absolute value equations} (AVE) and \textit{absolute value optimization} (AVO) problems
have been attracted much attention.
The AVE were introduced in 2004 by Rohn~\cite{rohn2004theorem}.
Basically, if $\tilde{A}, \tilde{B}$ are given matrices, and $\tilde{b}$  is a given vector,
one should find a vector $x$ that satisfies $\tilde{A} x + \tilde{B} |x| = \tilde{b}$,
where $|x|$ is a vector whose $i$-th entry is the absolute value of the $i$-th entry of $x$.
It is known that AVE are equivalent to the \emph{linear complementarity problems} (LCP)~\cite{hu2010note, mangasarian2006absolute, prokopyev2009equivalent},
which include many real-world applications.
As an extension of AVE, Mangasarian~\cite{mangasarian2007absolute} proposed in 2007
the AVO problems,
which have the absolute value of variables
in their objective and constraint functions.
More precisely, the AVO problem considered is given by
\[
\left.
\begin{array}{lrl}
		& \min	 	& \tilde{c}^T x + \tilde{d}^T |x| \\
 		& \mathrm{s.t.} & \tilde{A} x + \tilde{B} |x| = \tilde{b}, \\
		& 		& \tilde{H} x + \tilde{K} |x| \geq \tilde{p},
\end{array}
\right.
\]
where $\tilde{A}, \tilde{B}, \tilde{H}, \tilde{K}$ are given matrices,
and $\tilde{c}, \tilde{d}, \tilde{b}, \tilde{p}$  are vectors with appropriate dimensions.
Since AVE and LCP are equivalent, the AVO include
the \textit{mathematical programs with linear complementarity constraints}~\cite{luo1996mathematical},
which are one of the formulations of equilibrium problems.
As another application of AVO, Yamanaka and Fukushima~\cite{yamanaka2014branch}
presented facility location problems.

Since 2007, some methods for solving AVE have been presented in the literature.
For example, Rohn~\cite{rohn2009algorithm} considered an iterative algorithm
using the sign of variables for the case that $\tilde{A}$ and $\tilde{B}$ are square matrices.
For more general $\tilde{A}$ and $\tilde{B}$,
Mangasarian~\cite{mangasarian2007absolute} provided a method
involving successive linearization techniques.
Another methods include a concave minimization approach, given by Mangasarian~\cite{mangasarian2007absolute2},
and Newton-type methods, proposed by
Caccetta et al.~\cite{caccetta2011globally}, Mangasarian~\cite{mangasarian2009generalized},
and Zhang and Wei~\cite{zhang2009global}.
Some generalizations of AVE were also proposed.
For example, Hu et al.~\cite{hu2011generalized} considered an AVE involving
the absolute value of variables associated to the second-order cones.
Miao et al.~\cite{miao2015generalized} investigated an AVE with
the so-called circular cones.
In both papers, quasi-Newton based algorithms were used.

As for AVO problems, Yamanaka and Fukushima~\cite{yamanaka2014branch} proposed
to use a branch-and-bound technique.
In the branching procedure, two subproblems are generated by fixing the sign
of a variable as nonnegative or nonpositive.
In the bounding procedure, the dual information are considered.
However, to the best of our knowledge,
there is no other method that can find a global solution of AVO.
When comparing to AVE, the research associated to AVO problems is insufficient and
one of these reasons is the difficulty for obtaining feasible solutions of the problems.
In fact, their constraints include AVE, which are known to be NP-hard~\cite{mangasarian2007absolute}.

Another optimization problem that is related to AVO was recently investigated by
Friedlander et al.~\cite{FMP2014} and Aravkin et al.~\cite{ABD2017}.
It is called \emph{gauge optimization}, which basically consists in an optimization problem with the
so-called gauge function.
However, differently from AVO,
this problem does not consider multiple constraints, but only one gauge constraint.
In~\cite{ABD2017, FMP2014}, the authors showed that the Lagrange dual of gauge optimization problems
can be written in a closed-form by using the polar of the gauge functions.

In this paper, similarly to~\cite{ABD2017, FMP2014}, we introduce a generalized AVO problem, and show that it
has a wider practical application comparing to AVO problems.
It is also more general than gauge optimization problems,
because multiple constraints can be considered here.
The generalization is done by replacing absolute value functions
with \emph{positively homogeneous} functions.
So, the problem uses not only absolute value terms but also,
for instance, $p$-norm functions with $p \in (0, \infty ]$.
This generalized problem is referred here as \textit{positively homogeneous optimization}~(PHO).

Here, we introduce the PHO dual problem and compare it with the Lagrange dual.
We also show that the weak duality theorem holds,
similarly to the AVO problems~\cite{mangasarian2007absolute}.
In addition, we investigate the relation between the positively homogeneous duality
and the Lagrange duality, proving that these dual problems are equivalent
under some conditions.
In this case, the Lagrange dual of a positively homogeneous problem can be
written in a closed-form.
We point out that the gauge functions are special cases of the positively homogeneous functions,
which are not necessarily convex, differently from the gauge.
Moreover, the proposed problems here have
linear and positively homogeneous terms in their objective functions and constraints,
which is different from the problem considered in~\cite{ABD2017, FMP2014} that has only one gauge term.
Here, we also give some applications for the positively homogeneous problems,
which include $p$-order cone optimization,
sum of norms optimization and
group Lasso-type optimization problems,
and we show that their Lagrange dual can be written in a closed-form
even without convexity assumptions.

The paper is organized as follows. In Section 2, we give the definition of positively homoge\-neous functions
as well as its dual, showing some of their properties.
In Section 3, we define the PHO problems,
and we prove that weak duality holds.
In Section 4, the relation between the Lagrangian dual and the positively homogeneous dual is discussed.
We give some applications for PHO problems in Section~5.
We conclude the paper in Section 6, with final remarks and some future works.

We consider the following notations throughout the paper.
We denote by $\mb{R}_{++}$ the set of positive real numbers.
Let $x \in \mb{R}^n$ be a $n$-dimensional column vector, and $A \in \mb{R}^{n \times m}$ be a matrix
with dimension $n \times m$. We use $T$ to denote transpose.
For two vectors $x$ and $y$, we denote the vector $(x^T, y^T)^T$ as $(x, y)^T$ for simplicity.
If $x \in \mb{R}^n$, then its $i$-th entry is denoted by $x_i$, so $x=(x_1, \dots, x_n)^T$.
Moreover, if $I \subseteq \{1, \dots, n\}$, then $x_I$
corresponds to the subvector of $x$ with entries $x_i$, $i \in I$.
The notation $\# J$ denotes the number of elements of a set $J$.
The identity matrix with dimension $n$ is given by $E_n \in \mb{R}^{n \times n}$.
Also, we denote by $\| \cdot \|_p$ and $\| \cdot \|_\infty$
the $p$-norm with $p > 0$ and the supremum norm, respectively.
If no distinction is made for the norm, we just use the notation $\|\cdot \|$.


\section{Positively homogeneous functions}

In this section, we first introduce the definitions of positively homogeneous
and vector positively homogeneous functions.
Then, we define their dual,
which will be used to describe the dual of PHO problems.
Moreover, we show some properties associated to these functions.

\begin{defn}(Positively homogeneous functions)\label{def:nnsh}
A function $\psi \colon \mb{R}^n \rightarrow \mb{R}$ is
\emph{positively homogeneous} if the following inequality holds:
\[
\psi( \lambda x ) = \lambda \psi(x) \qquad \mathrm{for ~ all ~} x \in \mb{R}^n, \lambda \in \mb{R}_{++}.
\]
\end{defn}

\begin{defn}(Vector positively homogeneous functions)\label{def:vnnsh}
A mapping $\Psi \colon \mb{R}^n \rightarrow \mb{R}^m$ is a \emph{vector positively homogeneous} function
if the following property holds:
\[
\Psi(x)=
\left[
\begin{array}{c}
\psi_1(x_{I_1})  \\
\vdots  \\
\psi_m(x_{I_m}) 
\end{array}
\right] \qquad \mathrm{for ~ all ~} x \in \mb{R}^n,
\label{asm:psi2_eq}
\]
where $\psi_i \colon \mb{R}^{n_i} \rightarrow \mb{R}$ is a positively homogeneous function
for all $i= 1, \ldots ,m$, $n = n_1 + \cdots + n_m$,
$I_i \subseteq \{1, \ldots, n \}$ is a set of indices satisfying
\[
I_i \cap I_j = \emptyset, \quad i \neq j, \qquad \mathrm{and} \qquad \#I_i = n_i,
\]
and $x_{I_i} \in \mb{R}^{n_i}$ is a disjoint subvector of $x$.
\end{defn}

The above definition basically says that $\Psi$ is vector positively homogeneous
if its block components are all positively homogeneous.
We now introduce the dual function of $\psi $,
which can be seen as a generalization of the dual norm.
Similarly, we also define the dual of vector positively homogeneous functions.

\begin{defn}(Dual positively homogeneous functions)\label{def:df}
Let $\psi \colon \mb{R}^n \rightarrow \mb{R}$ be a positively homogeneous function.
Then, $\psi^* \colon \mb{R}^n \rightarrow \mb{R}\cup\{\infty\}$ defined by
\[
\psi^*(y) := \sup \{ x^T y \: | \: \psi(x) \leq 1 \} \qquad \mathrm{for ~ all ~} y \in \mb{R}^n
\]
 is called the \emph{dual positively homogeneous} function of $\psi$.
\end{defn}

Note that $\psi^*$ is convex from definition.
In fact, for all $y, z\in \mb{R}^n$ and $\alpha \in (0, 1)$, we have
\begin{eqnarray*}
\psi^*(\alpha y+ (1-\alpha) z) 	&=& \sup \{ x^T (\alpha y + ( 1 - \alpha ) z ) \: | \: \psi(x) \leq 1 \} \\
						&\leq& \alpha \sup \{x^T y \: | \: \psi(x) \leq 1\} + (1 - \alpha) \sup \{ x^Tz \: | \: \psi(x) \leq 1 \} \\
						&=& \alpha \psi^*(y)+(1 - \alpha) \psi^*(z).
\end{eqnarray*}

\begin{defn}(Dual vector positively homogeneous functions)\label{def:dvnnsh}
Let $\Psi \colon \mb{R}^n \rightarrow \mb{R}^m$ be a vector positively homogeneous function.
A function $\Psi^* \colon \mb{R}^n \rightarrow \mb{R}^m$ is
a \emph{dual vector positively homogeneous} function  associated to $\Psi$ if the following property holds:
\[
\Psi^*(y)=
\left[
\begin{array}{c}
\psi_1^*(y_{I_1})  \\
\vdots  \\
\psi_m^*(y_{I_m}) 
\end{array}
\right], \quad
i= 1, \ldots ,m, \qquad \mathrm{for ~ all ~} y \in \mb{R}^n
\label{asm:psi2_eq}
\]
where $\psi_i^* \colon \mb{R}^{n_i} \rightarrow \mb{R}$ is the dual of positively homogeneous function $\psi_i$
for each $i=1,\ldots m$.
\end{defn}

In this paper, we assume two conditions for positively homogeneous functions.

\begin{asm}\label{asm:nn}
Let $\Psi \colon \mb{R}^n \rightarrow \mb{R}^m$ be a vector positively homogeneous function as in Definition~\ref{def:vnnsh}.
Then, for all $i=1, \ldots, m$,
the positively homogeneous function $\psi_i$ satisfies the following conditions:
\begin{enumerate}
	\item	$\psi_i(x_{I_i}) \geq 0 \quad \mathrm{for ~ all ~} x_{I_i} \in \mb{R}^{n_i}$,
	\item	$\mathrm{If}$ $x_{I_i} \neq 0$, $\mathrm{then}$ $\psi_i(x_{I_i}) > 0$.
\end{enumerate}
\end{asm}

From the definition of positively homogeneous functions, we observe that $\psi_i(0) = 0$.
In fact, if $x = 0$ then $0 = \psi (\lambda x) - \lambda \psi(x) = (1 - \lambda) \psi_i(0)$ for all $\lambda \in \mb{R}_{++}$.
Moreover, the second condition of the above assumption shows that zero is the only point that satisfies $\psi_i(x) = 0$.
We also observe that if $\psi_i$ is taken as the usual vector norm, then it satisfies these assumptions.
Note that under the above assumption, the dual function $\psi_i^*$ always takes finite values.

We now show an important property satisfied by vector positively homogeneous functions
and their dual.

\begin{prop}\label{prop:ineq}
Let $\Psi$ and $\Psi^*$ be a vector positively homogeneous function and its dual, respectively.
Suppose that Assumption \ref{asm:nn} holds.
Then, the following inequalities hold:
\begin{eqnarray*}
\Psi^*(y) 	&\geq&  0,\\
\Psi(x)^T \Psi^*(y) &\geq&  x^T y
\end{eqnarray*}
for any $x, y \in \mb{R}^n$.
\end{prop}

\begin{prf}
For simplicity, we take an arbitrary index $i$ and denote $\psi_i$ and $x_{I_i}$ as $\psi$ and $x$, respectively.
From Definition \ref{def:nnsh}, we have $\psi(0) = 0$.
Using this result and Definition \ref{def:df}, we obtain
\[
\psi^*(y) = \sup \{ x^T y \: | \: \psi(x) \leq 1 \} \geq 0 \qquad \mathrm{for ~ all ~} y \in \mb{R}^n.
\]
This shows that $\Psi^*(y) \geq 0 ~ \mathrm{for ~ all ~} y \in \mb{R}^n$ from Definition \ref{def:dvnnsh}.

If $x = 0$, then the second inequality of this proposition clearly holds.
If $x \neq 0$, then $\psi(x) > 0$ from Assumption \ref{asm:nn} and so
\[
\psi \left( \dfrac{x}{\psi(x)} \right) = \dfrac{1}{\psi(x)} \psi(x) = 1
\]
holds once again from Definition \ref{def:nnsh}. Therefore, we obtain
\begin{eqnarray*}
\psi^*(y) \geq \left( \dfrac{x}{\psi(x)} \right)^T y \qquad \mathrm{for ~ all ~} y \in \mb{R}^n.
\end{eqnarray*}
Then, for all $x, y \in \mb{R}^n$, we have
\[
\psi(x) \psi^*(y) \geq x^T y,
\]
which indicates that
\[
\Psi(x)^T \Psi^*(y) = \sum_{i=1}^m \psi_{I_i}(x) \psi_{I_i}^*(y) \geq \sum_{i=1}^m x_{I_i}^T y_{I_i} = x^T y.
\]
\qed
\end{prf}


\section{Positively homogeneous optimization problems}

We consider the following \emph{positively homogeneous optimization} (PHO) problem:

\[
\left.
\begin{array}{lrl}
		& \min	 	& c^T x + d^T \Psi(x) \\
 		& \mathrm{s.t.} & A x + B \Psi(x) = b, \\
		& 		& H x + K \Psi(x) \geq p,
\end{array}
\right. \tag{$\mathrm{P}$}
\]
where $c\in \mb{R}^n, d\in \mb{R}^m, b\in \mb{R}^k, p\in \mb{R}^\ell,
A\in \mb{R}^{k\times n}, B\in \mb{R}^{k\times m}, H\in \mb{R}^{\ell \times n}$ and $K\in \mb{R}^{\ell \times m}$
are given constant vectors and matrices,
and $\Psi \colon \mb{R}^n \rightarrow \mb{R}^m$ is a vector positively homogeneous function
satisfying Assumption \ref{asm:nn}.

Now we give the Lagrangian dual of the problem (P) as follows:
\[
\sup_{\substack{u\\ v\geq 0}} \: \omega(u, v) , \tag{$\mathrm{D_{\cal L}}$}
\]
where $\omega \colon \mb{R}^k \times \mb{R}^\ell \rightarrow \mb{R}$ is given by
\begin{eqnarray}\label{omega}
\omega(u, v) := \inf_x {\cal L}(x, u, v),
\end{eqnarray}
and ${\cal L}\colon \mb{R}^n \times \mb{R}^k \times \mb{R}^\ell \rightarrow \mb{R}$ is the Lagrangian function of (P)
defined by
\begin{eqnarray*}
{\cal L}(x, u, v) &:=& c^T x + d^T \Psi(x) + u^T (b - A x - B \Psi(x)) + v^T (p - H x - K \Psi(x)) \\
 &=& b^T u + p^T v - (A^T u + H^T v - c)^T x + (d - B^T u - K^T v)^T \Psi(x),
\end{eqnarray*}
with $u\in \mb{R}^k$ and $v\in \mb{R}^\ell$ as the \textit{Lagrange multipliers}
associated to the equality and inequality constraints, respectively.
Notice that it is difficult to write concretely the objective function of the problem ($\mathrm{D_{\cal L}}$)
because it is, in general, not convex with respect to $x$.

In order to obtain a closed-form dual problem,
we consider a convex relaxation of the original problem (P) and its Lagrangian dual.
For simplicity, we investigate the case where $\Psi(x) = |x| := (|x_1|, \ldots, |x_n|)^T$, and
(P) has a linear objective function and only inequality constraints.
More precisely, we analyze the following problem:
\[
\left.
\begin{array}{lrl}
		& \min	 	& c^T x \\
 		& \mathrm{s.t.} & A x + B |x| \geq b.
\end{array}
\right. \tag{$\mathrm{P_{a}}$}
\]
If we set $x = x^+ - x^-$ and $|x| = x^+ + x^-$, where $x^+_i = \max \{0,  x_i \}$ and $x^-_i = \max \{0,  -x_i \}$,
then we can write ($\mathrm{P_{a}}$) as
\[
\left.
\begin{array}{lrl}
		& \min	 	& [c^T | -c^T]
\left[
\begin{array}{c}
x^+ \\
x^-
\end{array}
\right] \vspace{1mm} \\
 		& \mathrm{s.t.} & [A | -A]
\left[
\begin{array}{c}
x^+ \\
x^-
\end{array}
\right] + [B | B]
\left[
\begin{array}{c}
x^+ \\
x^-
\end{array}
\right] \geq b,
\end{array}
\right.
\]
which is equivalent to the following problem:
\[
\left.
\begin{array}{lrl}
		& \min	 	& [c^T | -c^T]
\left[
\begin{array}{c}
y_1 \\
y_2
\end{array}
\right] \vspace{1mm} \\
 		& \mathrm{s.t.} & [A | -A]
\left[
\begin{array}{c}
y_1 \\
y_2
\end{array}
\right] + [B | B]
\left[
\begin{array}{c}
y_1 \\
y_2
\end{array}
\right] \geq b, \\
		&			& y_1, y_2 \geq 0, \\
		&			& y_1^T y_2 = 0,
\end{array}
\right.
\]
where $y_1, y_2 \in \mb{R}^n$.
Notice that the above problem is not convex due to the complementarity constraint $y_1^T y_2 = 0$.
Therefore, we remove it from the problem and obtain the following relaxed one:
\[
\left.
\begin{array}{lrl}
		& \min	 	& [c^T | -c^T] ~y \\
 		& \mathrm{s.t.} & [A + B | - A + B] ~y \geq b, \\
		&			& y \geq 0,
\end{array}
\right.
\]
where $y = (y_1, y_2)^T$.
This problem is just a linear programming,
then its Lagrangian dual can be written easily as
\[
\left.
\begin{array}{lrl}
		& \max	 	& b^T u \\
 		& \mathrm{s.t.} & 
		\left[
		\begin{array}{c}
		A^T + B^T \\
		- A^T + B^T
		\end{array}
		\right] ~u \leq 
		\left[
		\begin{array}{c}
		c \\
		- c
		\end{array}
		\right], \\
		&			& u \geq 0.
\end{array}
\right.
\]
Observing that the first constraint is equivalent to $|A^T u - c | + B^T u \leq 0$,
we finally obtain the following closed-form dual problem:
\[
\left.
\begin{array}{lrl}
		& \max	 	& b^T u \\
 		& \mathrm{s.t.} & |A^T u - c | + B^T u \leq 0, \\
		&			& u \geq 0.
\end{array}
\right. \tag{$\mathrm{D_{a}}$}
\]
In fact, the problem ($\mathrm{D_{a}}$) is the AVO dual of ($\mathrm{P_{a}}$)
proposed by Mangasarian in \cite{mangasarian2007absolute},
and the weak duality clearly holds in this case.

Let us return to the general problem (P).
Inspired by the above AVO dual problem ($\mathrm{D_{a}}$),
we consider the following problem as the positively homogeneous dual problem:
\[
\left.
\begin{array}{lrl}
		& \max	 	& b^T u + p^T v \\
 		& \mathrm{s.t.} & \Psi^*(A^T u + H^T v - c) + B^T u + K^T v \leq d, \\
		& 		& v \geq 0,
\end{array}
\right. \tag{$\mathrm{D}$}
\]
where $\Psi^*$ is the dual vector positively homogeneous function associated to $\Psi$.
Note that (D) is a convex optimization problem
since each component $\psi^*_i$ of $\Psi^*$ is a convex function.

The theorem below shows that the proposed dual problem (D) is reasonable,
in the sense that the weak duality holds between (P) and (D).

\begin{thm}\label{prop:wd}(Weak duality)
For problems (P) and (D), the following inequality holds:
\[
c^T x + d^T \Psi(x) \geq b^T u + p^T v
\]
for all feasible points $x \in \mb{R}^n$ and $(u, v) \in \mb{R}^k \times \mb{R}^\ell$ of (P) and (D), respectively.
\end{thm}

\begin{prf}
Let $x \in \mb{R}^n$ and $(u, v) \in \mb{R}^k \times \mb{R}^\ell$ be feasible
for (P) and (D), respectively. Then, we have
\begin{eqnarray*}
c^T x + d^T \Psi(x) 	&\geq & c^T x + ( \Psi^*(A^T u + H^T v - c) + B^T u + K^T v )^T \Psi(x) \\
				&=	  & c^T x + \Psi^*(A^T u + H^T v - c)^T \Psi(x) + u^T B \Psi(x) + v^T K \Psi(x),
\end{eqnarray*}
where the inequality holds from the first constraint of (D) and the nonnegativity of $\Psi$.
From the second inequality of Proposition \ref{prop:ineq}, we also obtain:
\begin{eqnarray*}
c^T x + d^T \Psi(x)	&\geq	& c^T x + (A^T u + H^T v - c)^T x + u^T B \Psi(x) + v^T K \Psi(x) \\
				&=	& u^T (A x + B \Psi(x)) + v^T (H x + K \Psi(x)).
\end{eqnarray*}
Finally, the constraints of (P) gives
\[
c^T x + d^T \Psi(x) \geq b^T u + p^T v,
\]
which completes the proof.\qed
\end{prf}

The weak duality theorem itself is a powerful theoretical result,
but it does not mention how large the duality gap between (P) and (D) is.
And the duality gap can be large depending on problems,
then the dual problem (D) may be useless.
Therefore, in the next section,
we investigate the relation between the Lagrangian dual problem ($\mathrm{D_{\cal L}}$)
and the one (D) proposed here.
As a result, surprisingly, we find that ($\mathrm{D_{\cal L}}$) and (D) are equivalent.


\section{The positively homogeneous duality and the\\ Lagrangian duality}

In this section, we consider the relation between the positively homogeneous duality and
the more traditional Lagrangian duality of problem (P),
investigating conditions under which the Lagrangian dual problem $\mathrm{(D_{\cal L})}$
and the positively homogeneous dual problem (D) are equivalent.
Notice that the equivalence means the optimal values of (D) and $\mathrm{(D_{\cal L})}$
are the same if they are finite.
Recalling (\ref{omega}), we first show a condition that makes $\omega(\bar{u}, \bar{v})$,
the objective function of $\mathrm{(D_{\cal L})}$,
unbounded from below for some~$(\bar{u}, \bar{v})$.

\begin{lem}\label{lem:lag}
Let $\psi^*_i$ be the dual of the positively homogeneous functions $\psi_i$ for $i=1, \ldots, m$.
Suppose that Assumption \ref{asm:nn} holds. Also, assume that there exists $(\bar{u}, \bar{v})$ and an index $i_0$ satisfying
\[
\psi^*_{i_0} (\alpha_{I_{i_0}}) > \beta_{i_0},
\label{lem1}
\]
where $\alpha := A^T \bar{u} + H^T \bar{v} - c \in \mb{R}^n$, and $\beta := d - B^T \bar{u} - K^T \bar{v} \in \mb{R}^m$.
Then, there exists a sequence $\{ x^k \}$ such that $\|x^k\| \rightarrow +\infty$
and ${\cal L}(x^k, \bar{u}, \bar{v}) \rightarrow -\infty$ as $k \rightarrow +\infty$.
Therefore, $\omega(\bar{u}, \bar{v})$ is unbounded from below.
\end{lem}

\begin{prf}
Firstly, we denote $\bar{\alpha}$ and $\bar{\alpha}(\lambda)$ as follows:
\begin{eqnarray*}
\bar{\alpha} &:=& (\alpha_{I_1}, \alpha_{I_2}, \ldots, \alpha_{I_{i_0}}, \ldots, \alpha_{I_m}) \in \mb{R}^n, \\
\bar{\alpha}(\lambda) &:=& (\alpha_{I_1}, \alpha_{I_2}, \ldots, \lambda \hat{x}, \ldots, \alpha_{I_m}) \in \mb{R}^n,
\end{eqnarray*}
where $\lambda \in \mb{R}_{++}$ and $\hat{x} \in \mb{R}^{n_{i_0}}$ is defined as the supreme point of the following problem:
\[
\sup \{ x^T \alpha_{I_{i_0}} ~ | ~ \psi_{i_0}(x) \leq 1 \}.
\]
From the definition of $\hat{x}$, we obtain $\psi_{i_0} (\hat{x}) \leq 1$.
Then, from Definition \ref{def:df}, we have
\[
\hat{x}^T \alpha_{I_{i_0}} = \psi_{i_0}^*(\alpha_{I_{i_0}}) \geq \psi_{i_0}(\hat{x}) \psi_{i_0}^*(\alpha_{I_{i_0}}).
\]

The above equality and the definition of the Lagrangian function give
\begin{eqnarray*}
{\cal L}(\bar{\alpha}(\lambda), \bar{u}, \bar{v}) &=& b^T \bar{u} + p^T \bar{v} - \bar{\alpha}^T \bar{\alpha}(\lambda) + \beta^T \Psi(\bar{\alpha}(\lambda))  \\
	&= & b^T \bar{u} + p^T \bar{v} - \sum_{i\neq i_0} \alpha_{I_i}^T \alpha_{I_i} - \lambda \hat{x}^T \alpha_{I_{i_0}} + \sum_{i \neq i_0} \beta_i \psi_i(\alpha_{I_i}) + \beta_{i_0} \psi_{i_0}(\lambda \hat{x}) \\
        &= & \gamma - \lambda \hat{x}^T \alpha_{I_{i_0}} + \beta_{i_0} \psi_{i_0}(\lambda \hat{x}) \\
        &\leq & \gamma - \lambda \psi_{i_0}(\hat{x}) \psi_{i_0}^*(\alpha_{I_{i_0}}) + \beta_{i_0} \psi_{i_0}(\lambda \hat{x}),
\end{eqnarray*}
where $\gamma := b^T \bar{u} + p^T \bar{v} - \sum_{i\neq i_0} \alpha_{I_i}^T \alpha_{I_i} + \sum_{i \neq i_0} \beta_i \psi_i(\alpha_{I_i}) \in \mb{R}$ is constant with respect to $\lambda$.
Moreover, Definition \ref{def:nnsh} shows that
\begin{eqnarray*}
{\cal L}(\bar{\alpha}(\lambda), \bar{u}, \bar{v}) &=& \gamma - \lambda \psi_{i_0}(\hat{x}) \psi_{i_0}^*(\alpha_{I_{i_0}}) + \lambda \beta_{i_0} \psi_{i_0}(\hat{x}) \\
	&=& \gamma + \lambda \psi_{i_0}(\hat{x}) ( \beta_{i_0} - \psi^*_{i_0}(\alpha_{I_{i_0}})) \\
	&\leq& \gamma + \lambda ( \beta_{i_0} - \psi^*_{i_0}(\alpha_{I_{i_0}})).
\end{eqnarray*}
Therefore, ${\cal L}(\bar{\alpha}(\lambda), \bar{u}, \bar{v})$
converges to minus infinity when $\lambda $ increases.
Finally, if we set $x^k = \bar{\alpha}(\lambda^k)$ where $\lambda^k \rightarrow +\infty$ as $k \rightarrow +\infty$,
then ${\cal L}(x^k, \bar{u}, \bar{v}) \rightarrow -\infty$ and
we complete the proof.
\qed
\end{prf}

We now show that the positively homogeneous dual problem (D) and the Lagrangian one 
$\mathrm{(D_{\cal L})}$ are equivalent under some conditions.

\begin{lem}\label{thm:suf}
Suppose that Assumption \ref{asm:nn} holds. Assume also that the positively homogeneous dual problem (D) has a feasible solution
$(\bar{u}, \bar{v})  \in \mb{R}^k \times \mb{R}^\ell$,
and that there exists $x^* \in \mb{R}^n$ satisfying the following equality:
\begin{eqnarray}
(d - B^T \bar{u} - K^T \bar{v})^T \Psi(x^*) - (A^T \bar{u} + H^T \bar{v} - c)^T x^* = 0. \label{eq:suf}
\end{eqnarray}
Then, the positively homogeneous dual problem (D)
and the Lagrangian dual problem $(D_{\cal L})$ are equivalent.
\end{lem}

\begin{prf}
From Lemma \ref{lem:lag}, the function $\omega$ is unbounded from below
if there exists an index $i_0$ such that $\psi^*_{i_0} (\alpha_{I_{i_0}}) > \beta_{i_0}$,
where $\alpha := A^T \bar{u} + H^T \bar{v} - c \in \mb{R}^n$, and $\beta := d - B^T \bar{u} - K^T \bar{v} \in \mb{R}^m$.
Therefore, the problem $\mathrm{(D_{\cal L})}$ is equivalent to
\[
\left.
\begin{array}{lrl}
		& \sup	 	& \omega(u, v) \\
 		& \mathrm{s.t.} & \Psi^*(A^T u + H^T v - c) \leq d - B^T u - K^T v, \\
		& 		& v \geq 0.
\end{array}
\right. \tag{$\mathrm{D_{\cal L}'}$}
\]
Let $(\bar{u}, \bar{v})  \in \mb{R}^k \times \mb{R}^\ell$ be the feasible solution of $\mathrm{(D_{\cal L}')}$.
From the definition of the Lagrangian function, we obtain:
\begin{eqnarray*}
{\cal L}(x, \bar{u}, \bar{v}) &=& c^T x + d^T \Psi( x) + \bar{u}^T (b - A x - B \Psi( x )) + \bar{v}^T (p - H x - K \Psi( x )) \\
 &=& b^T \bar{u} + p^T \bar{v} - (A^T \bar{u} + H^T \bar{v} - c)^T x + (d - B^T \bar{u} - K^T \bar{v})^T \Psi( x ).
\end{eqnarray*}
Then, taking $x^* \in \mb{R}^n$ that satisfies (\ref{eq:suf}), we have
\[
{\cal L}(x^*, \bar{u}, \bar{v}) = b^T \bar{u} + p^T \bar{v}.
\]
Notice that $x^*$ is the solution of the problem
\[
\inf_x {\cal L}(x, \bar{u}, \bar{v}),
\]
because ${\cal L}(x, \bar{u}, \bar{v}) \geq b^T \bar{u} + p^T \bar{v} $
holds from Proposition \ref{prop:ineq}.
Therefore, the problem $\mathrm{(D_{\cal L}')}$ can be described as follows:
\[
\left.
\begin{array}{lrl}
		& \sup	 	& b^T u + p^T v \\
 		& \mathrm{s.t.} & \Psi^*(A^T u + H^T v - c) \leq d - B^T u - K^T v, \\
		& 		& v \geq 0,
\end{array}
\right.
\]
which is equivalent to the positively homogeneous dual problem $\mathrm{(D)}$.
\qed
\end{prf}

As a consequence of the above lemma,
we obtain the following result.

\begin{thm}\label{thm:sum}
Suppose that the Lagrangian dual problem ($D_{\cal L}$) has a feasible solution.
Assume also that the vector positively homogeneous function $\Psi$ satisfies Assumption \ref{asm:nn}.
Then, the positively homogeneous dual problem (D) and the Lagrangian dual problem ($D_{\cal L}$)
have the same optimal value and solutions.
\end{thm}

\begin{prf}
From Definition \ref{def:nnsh} and Assumption \ref{asm:nn}, we have $\Psi(0) = 0$.
It means that equation (\ref{eq:suf}) holds at $x^* = 0$.
Thus, from Lemma \ref{thm:suf}, the problems (D) and $\mathrm{(D_{\cal L})}$ have the same optimal value.

Moreover, we denote $S_D$ and $S_{D_{\cal L}}$ as the sets of
optimal solutions of problems (D) and ($\mathrm{D_{\cal L}}$), respectively.
Let us take $(u^*, v^*) \in S_D$. Then, it is clearly feasible for ($\mathrm{D_{\cal L}}$).
It follows from Theorem \ref{thm:sum} that the optimal values of (D) and ($\mathrm{D_{\cal L}}$) are the same,
which is $b^T u^* + p^T v^*$, and so $(u^*, v^*) \in S_{D_{\cal L}}$.
Conversely, let us take $(\bar u, \bar v) \in S_{D_{\cal L}}$.
Then, the point $(\bar u, \bar v)$ is feasible for ($\mathrm{D_{\cal L}}$).
Note that Lemma~\ref{lem:lag} indicates that if $(u, v)$ is feasible for ($\mathrm{D_{\cal L}}$)
and the objective function value of ($\mathrm{D_{\cal L}}$) at the point $(u, v)$ is finite,
then it is also feasible for (D).
Thus, $(\bar u, \bar v)$ is feasible for (D). Once again from Theorem \ref{thm:sum},
the optimal values of (D) and ($\mathrm{D_{\cal L}}$) are the same, which means that $(\bar u, \bar v) \in S_D$.
Consequently, we obtain $S_D = S_{D_{\cal L}}$.
\qed
\end{prf}

The above theorem shows that
the Lagrangian dual problem $\mathrm{(D_{\cal L})}$ can be written in a closed-form
when the function $\Psi$ is positively homogeneous and satisfies Assumption \ref{asm:nn}.
The paper \cite{mangasarian2007absolute} does not show that the same property
holds for the AVO problem. We now give it as a direct consequence of Theorem \ref{thm:sum}.

\begin{cor}\label{cor:avo}
If the dual of an AVO problem has a feasible solution,
then it is equivalent to the Lagrangian dual problem ($D_{\cal L}$).
\end{cor}

\begin{prf}
It holds from Theorem \ref{thm:sum} and the fact that the absolute value function is
positively homogeneous and satisfies Assumption \ref{asm:nn}.
\qed
\end{prf}

\begin{cor}
If the optimal values of an AVO primal problem and its Lagrangian dual problem ($D_{\cal L}$) are the same,
then the strong duality holds between the AVO primal and the AVO dual problem.
\end{cor}

\begin{prf}
It holds straightforward from Corollary \ref{cor:avo}.
\qed
\end{prf}

From the above result, AVO can be applied to
solve 0-1 integer optimization problems.
To solve such problems, their Lagrangian dual are often considered,
which is, in general, nondifferentiable due to the integer constraints.
On the other hand,
a 0-1 integer constraint, that is $x \in \{0, 1\}$, is equivalent to
$|2 x - 1 | = 1$. Then, 0-1 integer optimization problems
can be reduced to AVO,
and we obtain their AVO dual, which are actually \emph{linear programming} (LP) problems.
These LP dual problems are much easier to solve compared to the nondifferentiable ones. 
Therefore, it might be worth considering AVO dual problems from
the computational point of view.

\section{Examples of positively homogeneous optimization\\ problems}

In this section, we present several applications that are formulated as PHO,
and show their closed-form dual problems.

First, we observe that any $p$-norm function with $p \in [1, \infty)$ is positively homogeneous.
So, if $\psi$ is the $p$-norm, then $\psi^*$ becomes the $q$-norm, where $1/p + 1/q = 1$.
Therefore, if $\psi$ is taken as $\|\cdot \|_1, \|\cdot \|_2, \|\cdot \|_\infty$,
then $\psi^*$ becomes $\|\cdot \|_\infty, \|\cdot \|_2, \|\cdot \|_1$, respectively.
Moreover, in the case that $p \in (0, 1)$, the dual function $\psi^*$ is equal to $\| \cdot \|_\infty$ for all $p \in (0, 1)$,
which is proved in Proposition 6 of Appendix A.
From the result, we can consider any $p$-norm functions as $\psi$ in PHO problems.
And, even if such functions are nonconvex with $p \in (0, 1)$,
the Lagrangian dual problem can be written in a closed-form
from Theorem \ref{thm:sum}.

We now show some positively homogeneous problems using these $p$-norm functions.
The first example is the so-called linear second-order cone optimization problem \cite{alizadeh2003second},
which is one of the famous convex optimization problem.

\begin{exm}\label{exm:1}
Let $x=(x_1, x_2)^T \in \mb{R} \times \mb{R}^{n-1}$.
Then, we consider the linear second-order cone optimization problem written by
\[
\left.
\begin{array}{lrl}
				& \min	 	& c^T x \\
 				& \mathrm{s.t.} & A x = b, \\
				& 		& x_1 - \|x_2 \|_2 \geq 0,
\end{array}
\right. \tag{$\mathrm{P_1}$}
\]
where $c \in \mb{R}^n$, $A \in \mb{R}^{m \times n}$ and $b \in \mb{R}^m$.
The above problem can be written in PHO form as
\[
\left.
\begin{array}{lrl}
				& \min	 	& c^T x + 0^T \Psi(x) \\
 				& \mathrm{s.t.} & A x + 0 \Psi(x) = b, \\
				& 		&	H x + K \Psi(x) \geq 0,
\end{array}
\right.
\]
with $H = (1, 0, \ldots, 0) \in \mb{R}^{1 \times n}, K = (0, -1) \in \mb{R}^{1 \times 2}$ and $\Psi \colon \mb{R}^n \rightarrow \mb{R}^2, \Psi(x) = (|x_1|, \|x_2\|_2)^T$.
Then, recalling (D), its dual problem is given by
\[
\left.
\begin{array}{lrl}
				& \max	 	& b^T u \\
 				& \mathrm{s.t.} & \Psi^*(A^T u + H^T v - c) + K^T v \leq 0, \\
				& 			& v \geq 0,
\end{array}
\right.
\]
where $\Psi^* $ is identical to $\Psi$ in this case. Then, from the definition of $\Psi$, we have
\[
\left.
\begin{array}{lrl}
				& \max	 	& b^T u \\
 				& \mathrm{s.t.} & |(A^T u)_1 + v - c_1 | \leq 0, \\
				&			& \|(A^T u)_2 - c_2  \|_2 \leq v, \\
				& 			& v \geq 0,
\end{array}
\right.
\]
with $(A^T u)_1$ as the first component of $A^T u$, $(A^T u)_2$ is the rest of it,
and $c = (c_1, c_2)^T \in \mb{R} \times \mb{R}^{n-1}$.
The first constraint of the above problem shows that
\[
v = c_1 - (A^T u)_1,
\]
and $v \geq 0$ automatically holds from the second constraint.
Then, we obtain
\[
\left.
\begin{array}{lrl}
				& \max	 	& b^T u \\
 				& \mathrm{s.t.} & \|(A^T u)_2 - c_2  \|_2 \leq c_1 - (A^T u)_1
\end{array}
\right.
\]
as the dual problem of $\mathrm{(P_1)}$.
In fact, the above problem is the standard dual of the linear second-order cone optimization problem \cite{alizadeh2003second}.
\end{exm}

Although we use the $2$-norm in the above example,
any $p$-norm function with $p \in (0, \infty ]$ can be considered.
In this case, if $p \in [1, \infty ]$, then the primal and dual problems are
$p$-order cone and $q$-order cone optimization problems, respectively, where $1/p + 1/q = 1$ \cite{xue2000efficient}.
If $p \in (0, 1)$, then the dual is $\infty $-order cone optimization problem.

In the next example, we consider a gauge optimization problem,
which is also a convex problem with multiple gauge functions in its objective and constraint functions.
Here, we recall that $f$ is a gauge function if and only if
it is nonnegative, convex, positively homogeneous and satisfies $f(0)=0$ \cite{F87}.
For such a problem, we introduce its dual in PHO form.

\begin{exm}\label{exm:2}
Let $x \in \mb{R}^n$. We consider the following problem:
\[
\left.
\begin{array}{rl}
\min	 	& \displaystyle{\sum_{i=1}^s \alpha_i f_i (A_i x - a_i)} \\
\mathrm{s.t.} & g_j (B_j x - b_j) \leq \beta_j, \quad j=1,\ldots, t,
\end{array}
\right.
\tag{$\mathrm{P_2}$}
\]
where $\alpha_i, \beta_j \in \mb{R}_+$, $A_i \in \mb{R}^{m_i \times n}$, $B_j \in \mb{R}^{k_j \times n}$,
$a_i \in \mb{R}^{m_i}$ and $b_j \in \mb{R}^{k_j}$ are given for all $i=1, \ldots, s$ and $j=1, \ldots, t$,
and $f_i \colon \mb{R}^{m_i} \rightarrow \mb{R}$ and $g_j \colon \mb{R}^{k_j} \rightarrow \mb{R}$ are gauge functions.
Letting $y_i := A_i x - a_i$ and $z_j := B_j x - b_j$, $\mathrm{(P_2)}$ can be written as
\[
\left.
\begin{array}{rll}
\min	 	& \displaystyle{\sum_{i=1}^s \alpha_i f_i (y_i)} & \\
\mathrm{s.t.} & g_j (z_j) \leq \beta_j, 	& j=1,\ldots, t,\\
			& A_i x - y_i = a_i, 	& i=1,\ldots, s,\\
			& B_j x - z_j = b_j, 	& j=1,\ldots, t.
\end{array}
\right.
\]
The above problem does not have a gauge function defined for the variable $x$,
so we introduce such a gauge function $x \mapsto \psi(x)$ and
rewrite the problem into the following way:
\[
\left.
\begin{array}{rl}
\min	 	& 0 \times \psi(x) + \displaystyle{\sum_{i=1}^s \alpha_i f_i (y_i)} + 0 \times \displaystyle{\sum_{j=1}^t g_j (z_j)} \\
\mathrm{s.t.}	& 0 \times \psi(x) \leq 0, \\
			& 0 \times f_i(y_i) \leq 0,	\hspace{5mm}	i=1,\ldots, s, \\
			& g_j (z_j) \leq \beta_j, 	\hspace{10mm}	j=1,\ldots, t, \\
			& A_i x - y_i = a_i, 		\hspace{6mm}	i=1,\ldots, s,\\
			& B_j x - z_j = b_j, 		\hspace{5mm}	j=1,\ldots, t.
\end{array}
\right.
\]
Note that $\psi \colon \mb{R}^{n} \rightarrow \mb{R}$ is a dummy gauge function with $x$ as its domain.

Let
\[
\hat{x} := (x, y_1, \ldots, y_s, z_1, \dots, z_t) \in \mb{R}^{n + \sum_{i=1}^s m_i + \sum_{j=1}^t k_j}
\]
and
\[
\Psi(\hat{x}) := (\psi(x), f_1(y_1), \ldots, f_s(y_s), g_1(z_1), \ldots, g_t(z_t))^T.
\]
Then the above problem can be rewritten as
\[
\left.
\begin{array}{rl}
\min	 	& d^T \Psi(\hat{x}) \\
\mathrm{s.t.} & K \Psi(\hat{x}) \leq p, \\
			& \hat{A} \hat{x} = \hat{b},
\end{array}
\right.
\]
where $d = (0, \alpha_1, \ldots, \alpha_s, 0, \ldots, 0)^T \in \mb{R}^{1 + s + t}$,
$p = (0, \ldots, 0, \beta_1, \ldots, \beta_t)^T \in \mb{R}^{1 + s + t}$,
\[
K = 
\left[
\begin{array}{cc}
 0		& 0	 \\
 0		& E_t	
\end{array}
\right],
\hat{A} = 
\left[
\begin{array}{c|ccc|ccc}
A_1		& - E_{m_1} 	& 		&			&		&		&	 \\
\vdots	& 		 	& \ddots	&			&		& 0		&	\\
A_s		&			&		& -E_{m_s} 	&		&		&	\\ \hline
B_1		&			&		&			& -E_{k_1}&		&	\\
\vdots	&			& 0		&			&		& \ddots	&	\\
B_t		&			&		&			&		&	 	&	-E_{k_t}
\end{array}
\right], and \:
\hat{b} = 
\left[
\begin{array}{c}
a_1	 \\
\vdots \\
a_s	\\ \hline
b_1	\\
\vdots	\\
b_t	
\end{array}
\right].
\]
Moreover, its positively homogeneous dual problem is given by
\[
\left.
\begin{array}{rl}
\max	 	& \hat{b}^T u - p^T v \\
\mathrm{s.t.} & \Psi^*(\hat{A}^T u) - K^T v \leq d, \\
		& v \geq 0.
\end{array}
\right.
\]
For simplification, let $u = (u_{11}, \ldots, u_{1s}, u_{21}, \ldots, u_{2t})^T$
with $u_{1i} \in \mb{R}^{m_i}, i=1, \ldots, s$
and $u_{2j} \in \mb{R}^{k_j}, j=1, \ldots, t$.
Then the above problem is rewritten as
\[
\left.
\begin{array}{rl}
\max	 	& \displaystyle{ \sum_{i=1}^s a_i^T u_{1i} + \sum_{j=1}^t b_j^T u_{2j} - \sum_{\ell=1}^{t} \beta_\ell v_{1+s+\ell} } \\
\mathrm{s.t.} 	& \displaystyle{\sum_{i=1}^s A_i^T u_{1i} + \sum_{j=1}^t B_j^T u_{2j} = 0}, \\
			& f_i^* (-u_{1i}) \leq \alpha_i, 	\hspace{11mm} i=1,\ldots, s, \\
			& g_j^* (-u_{2j}) \leq v_{1+s+j}, \quad j=1,\ldots, t.
\end{array}
\right.
\tag{$\mathrm{D_{2}}$}
\]
Notice that the last constraint implies $v \geq 0$
because $g_j^*$ is also a gauge function.
Moreover, $\mathrm{(D_2)}$ does not include 
the dual function $\psi^*$ of the dummy gauge function $\psi$.
\end{exm}

The next example is the group Lasso-type problems \cite{meier2008group,yuan2006model},
which is a special case of $\mathrm{(P_2)}$
and consist in unconstrained minimizations of the sum of certain norms.
Such problems have many applications, in particular
they appear in compressed sensing area \cite{eldar2009robust,stojnic2009reconstruction},
where the sparsity of solutions are important.
As an example, we consider a primal problem with $p_1$-norm and $p_2$-norm where $p_1, p_2 \in \mb{R}_+$,
which are used in the regularization terms.

\begin{exm}\label{exm:3}
Let $x \in \mb{R}^n$ and $p_1, p_2 \in \mb{R}_+$. We consider the following problem:
\[
\min \quad  \| A x - b \|_2 + \lambda_1 \sum_{i=1}^{m'} \| x_{I_i} \|_{p_1} + \lambda_2 \sum_{i=m'+1}^{m} \| x_{I_i} \|_{p_2}
\tag{$\mathrm{P_3}$}
\]
where $\lambda_1, \lambda_2 \in \mb{R}_+, b \in \mb{R}^m, A \in \mb{R}^{m \times n}$ and $0 < m' < m$.

Notice that the first term of the objective function of group Lasso-type problems
are usually the square of $2$-norm functions.
However, it is not positively homogeneous,
so we removed the square and considered just the $2$-norm functions.

We obtain the above problem by setting, in $\mathrm{(P_2)}$, $s = m + 1$,
\[
\alpha_i =
\left\{
\begin{array}{ll}
\lambda_1, &	\: \mathrm{if} \quad i = 1, \ldots, m',	\\
\lambda_2, &	\: \mathrm{if} \quad i = m'+1, \ldots, m,	\\
1,	 	 &	\: \mathrm{if} \quad i = m+1,
\end{array}
\right.
\]
\[
A_i =
\left\{
\begin{array}{ll}
E_{I_i}, &	\: \mathrm{if} \quad i = 1, \ldots, m, \\
A,	 	 &	\: \mathrm{if} \quad i = m + 1,
\end{array}
\right.
\]
where $E_{I_i}$ is a submatrix of $E_n$ with $E_j,~ j\in I_i$ as its rows,
\[
a_i =
\left\{
\begin{array}{ll}
0, &	\: \mathrm{if} \quad i = 1, \ldots, m, \\
b,	 	 &	\: \mathrm{if} \quad i = m+1,
\end{array}
\right.
\]
and
\[
f_i(\cdot) =
\left\{
\begin{array}{ll}
\| \cdot \|_{p_1}, &	\: \mathrm{if} \quad i = 1, \ldots, m',	\\
\| \cdot \|_{p_2}, &	\: \mathrm{if} \quad i = m'+1, \ldots, m,	\\
\| \cdot \|_2, 	 &	\: \mathrm{if} \quad i = m+1.
\end{array}
\right.
\]
Then, recalling $\mathrm{(P_2)}$ and $\mathrm{(D_2)}$, the dual of $\mathrm{(P_3)}$ can be written as
\[
\left.
\begin{array}{rl}
\max	 	& b^T u_{1(m+1)} \\
\mathrm{s.t.} & \displaystyle{ \sum_{i=1}^{m} E_{I_i}^T u_{1i} + A^T u_{1(m+1)} = 0 }, \\
		& \| -u_{1i} \|_{q_1} \leq \lambda_1, \quad i = 1,\dots ,m', \\
		& \| -u_{1i} \|_{q_2} \leq \lambda_2, \quad i = m'+1,\dots ,m, \\
		& \| -u_{1(m+1)} \|_{2} \leq 1,
\end{array}
\right.
\]
where $q_i, i=1, 2$ are obtained by
\[
q_i = 
\left\{
\begin{array}{ll}
\dfrac{p_i}{p_i - 1},	& \: \mathrm{if} \quad p_i > 1, 	 \\
\infty ,			& \: \mathrm{if} \quad p_i \in (0, 1],
\end{array}
\right. \tag{3}
\]
from Proposition \ref{prop:exf2} of Appendix A.
Notice that the first equality constraint can be rewritten~as
\[
u_{1i} + (A^T)_{I_i} u_{1(m+1)} = 0, \quad i=1, \ldots, m.
\]
Then, the above problem is described as
\[
\left.
\begin{array}{rl}
\max	 	& b^T u \\
\mathrm{s.t.} &\| (A^T)_{I_i} u \|_{q_1} \leq \lambda_1, \quad i = 1,\dots ,m', \\
		& \| (A^T)_{I_i} u \|_{q_2} \leq \lambda_2, \quad i = m'+1,\dots ,m, \\
		&  \| -u \|_{2} \leq 1,
\end{array}
\right.
\]
where we denote $u_{1(m+1)}$ as $u$ for simplicity.
\end{exm}

The next example is also a Lasso-type problem.
In this case, the objective function is a gauge, because the sum of gauge functions is also gauge.
In order to obtain the dual of a gauge optimization problem,
the polar of the objective function should be considered \cite{ABD2017, FMP2014}.
However, it may be difficult to obtain the polar of a sum of gauge functions.
To overcome this drawback, we use here the PHO framework.

\begin{exm}\label{exm:4}
Let $x \in \mb{R}^n$ and $p_1, p_2 \in \mb{R}_+$. We consider the following problem:
\[
\left.
\begin{array}{rl}
\min	 	& \lambda_1 \| x \|_{p_1} + \lambda_2 \| x \|_{p_2} \\
\mathrm{s.t.} 	& \| A x - b \|_2 \leq \beta,
\end{array}
\right.
\tag{$\mathrm{P_4}$}
\]
where $\lambda_1, \lambda_2, \beta \in \mb{R}_+$, $A \in \mb{R}^{m \times n}$ and $b \in \mb{R}^m$.
The above problem can be obtained if we set, in $\mathrm{(P_2)}$, $s = 2$, $t=1$,
$\alpha_1 = \lambda_1$, $\alpha_2 = \lambda_2$,
$A_1 = A_2 = E_n$, $a_1 = a_2 = 0$, $B_1 = A$, $b_1 = b$,
$f_1(\cdot) = \| \cdot \|_{p_1}$, $f_2(\cdot) = \| \cdot \|_{p_2}$, $g_1(\cdot) = \| \cdot \|_2$.
Then, recalling $\mathrm{(D_2)}$, the dual of $\mathrm{(P_4)}$ is written by
\[
\left.
\begin{array}{rl}
\max	 	& b^T u_{21} - \beta v_4 \\
\mathrm{s.t.} 	& u_{11} + u_{12} + A^T u_{21} = 0, \\
			& \|-u_{11} \|_{q1} \leq \lambda_1, \\
			& \|-u_{12} \|_{q2} \leq \lambda_2, \\
			& \|-u_{21} \| \leq v_4,
\end{array}
\right.
\]
which is finally rewritten as
\[
\left.
\begin{array}{rl}
\max	 	& b^T u_{2} - \beta v \\
\mathrm{s.t.} 	& \| u_{1} + A^T u_{2} \|_{q1} \leq \lambda_1, \\
			& \|-u_{1} \|_{q2} \leq \lambda_2, \\
			& \|-u_{2} \| \leq v,
\end{array}
\right.
\]
where we set $u_{12}, u_{21}$ and $v_4$ as $u_1, u_2$ and $v$, respectively,
and $q_1$ and $q_2$ are defined in $\mathrm{(3)}$.

\end{exm}

In order to control the sparsity of the solutions of the above Lasso-type problems,
we can use any combination of $p$-norm functions, with $p \in (0, \infty ]$,
as the regularization terms.
Especially, it is reported that the $p$-norm functions with $p \in (0, 1)$ in $\mathrm{(P_3)}$ is useful
because they give sparser solutions than $1$-norm functions \cite{C2007, CY2008, MR2010}.

We now give another example: the sum of norms optimization problems,
which are generally nonconvex.
Such problems have applications, for example, in facility location,
where locations of new facilities should be decided by analyzing
the distance between the new and the existing facilities \cite{w2011}.
Moreover, the problem of the following example can be applied
not only to the minimization of the distance but also maximization of it
by taking the constant $\lambda_i$ as $-\lambda_i$.
Such a situation can be found for instance in locating obnoxious facilities in residential areas.

\begin{exm}\label{exm:5}
Let $x \in \mb{R}^n$. We consider the following problem:
\[
\left.
\begin{array}{rl}
\min	 	& \displaystyle{ \sum_{i=1}^{s} \lambda_i f_i( A_i x - a_i), } \\
\mathrm{s.t.} & B x \leq b,
\end{array}
\right.
\tag{$\mathrm{P_5}$}
\]
where $\lambda_i \in \mb{R}$, $A_i \in \mb{R}^{m_i \times n}$, $B \in \mb{R}^{k \times n}$,
$a_i \in \mb{R}^{m_i}$ and $b \in \mb{R}^{k}$ are given,
and $f_i \colon \mb{R}^{m_i} \rightarrow \mb{R}$, $i=1, \ldots, s$ are positively homogeneous functions.
We now introduce its positively homogeneous dual
by taking almost the same procedure as in Example 2.
Let $y_i := A_i x - a_i$, then $\mathrm{(P_5)}$ is equivalent to
\[
\left.
\begin{array}{rl}
\min	 	& \displaystyle{\sum_{i=1}^{s} \lambda_i f_i ( y_i ) } \\
\mathrm{s.t.} & A_i x - y_i = a_i, \quad i=1,\ldots s, \\
			& B x \leq b.
\end{array}
\right.
\]
By introducing additional constraints, we consider the following problem:
\[
\left.
\begin{array}{rl}
\min	 	& \displaystyle{\sum_{i=1}^{s} \lambda_i f_i ( y_i ) } \\
\mathrm{s.t.} & A_i x - y_i = a_i, \quad i=1,\ldots s, \\
			& B x \leq b, \\
			& c_i f_i(y_i) \leq d_i, \hspace{7mm} i=1, \ldots s,
\end{array}
\right.
\]
where $c_i$ and $d_i$ are strictly positive constants.
Notice that the additional constraints ensure the boundedness of
the each term of the objective function
especially when $\lambda_i$ is strictly negative.
Without such constraints,
$(\mathrm{P_5})$ can be unbounded depending on the linear constraint,
and then its dual becomes infeasible.
Note that the additional constraints do not change solutions,
when we choose $c_i$ and $d_i$ so that the constraint $c_i f_i(y_i) \leq d_i$
will include reasonable solutions.

Let $\hat{x} := (x, y_1, \ldots, y_s)^T \in \mb{R}^{n + \sum_{i=1}^s m_i}$
and $\Psi(\hat{x}) := (\psi(x), f_1(y_1), \ldots, f_s(y_s))^T \in~\mb{R}^{1 + s}$,
where $\psi(\cdot)$ is a dummy positively homogeneous function.
Then the above problem can be described as
\[
\left.
\begin{array}{rl}
\min	 	& d^T \Psi(\hat{x}) \\
\mathrm{s.t.} & \hat{A} \hat{x} = \hat{a}, \\
			& H \hat{x} + K \Psi(\hat{x}) \geq p,
\end{array}
\right.
\]
where $d = (0, \lambda_1, \ldots, \lambda_s)^T$, $\hat{a} = (a_1, \ldots, a_s)^T$, $p = (-b, -d_1, \ldots, -d_s)^T$,
\[
\hat{A} = \left[
    \begin{array}{cccc}
      A_1 	& -E_{m_1} 	& 		&0  \\
      \vdots &   		& \ddots 	&  \\
      A_s 	&   0		& 		& -E_{m_s}
    \end{array}
  \right], \:
H = 
\left[
\begin{array}{cc}
-B		& 0		 \\
0		& 0			
\end{array}
\right], \: and \:
K = 
\left[
\begin{array}{c|ccc}
0		&  			& 	0	&			 \\ \hline
		& 	-c_1	 	& 		&			\\
0		&			&	\ddots	& 		\\
		&			&			& -c_s		
\end{array}
\right].
\]
Then, recalling the positively homogeneous dual (D),
the dual of the above problem can be written~as
\[
\left.
\begin{array}{rl}
\max	 	& \hat{a}^T u + p^T v \\
\mathrm{s.t.} & \Psi^*(\hat{A}^T u + H^T v) \leq d - K^T v, \\
			& v \geq 0,
\end{array}
\right.
\]
which is rewritten by
\[
\left.
\begin{array}{rl}
\max	 	& \displaystyle{ \sum_{i=1}^s a_i^T u_i - b^T v_1 - \sum_{i=1}^s d_i^T v_{i+1}} \\
\mathrm{s.t.} & \displaystyle{\sum_{i=1}^s A_i^T u_i - B^T v_1 = 0}, \\
			& f_i^* ( -u_i ) \leq \lambda_i + c_i, , \quad i=1, \ldots s, \\
			& v \geq 0,
\end{array}
\right.
\]
where $v=(v_1, \ldots, v_{s+1})^T$.
\end{exm}


\section{Conclusion}

In this paper, we proposed an optimization problem with positively homogeneous functions,
which we call positively homogeneous optimization problem.
We also introduced its dual problem and showed the weak duality theorem between these problems.
Moreover, we gave sufficient conditions for the equivalency between the proposed dual
and the Lagrangian dual problems.
Finally, we presented some examples of positively homogeneous problems
to show their value in real-world applications.
One natural future work will be to propose methods that obtain approximate solutions
of positively homogeneous optimization problems.
We believe the theoretical results described here are essential for that.

\vspace{5mm}
\noindent
\textbf{Acknowledgements}
The authors are grateful to Prof. Ellen. H. Fukuda for helpful comments and suggestions.
This work was supported in part by a Grant-in-Aid for Scientific Research (C) (17K00032) from 
Japan Society for the Promotion of Science.


\appendix

\section{Appendix}

The following proposition shows that
the dual of the $p$-norm function is the $\infty$-norm
even when $p$ is less than 1.

\begin{prop}\label{prop:exf2}
Suppose that $p \in (0, 1)$.
Then, the dual of the $p$-norm function is equal to the $\infty$-norm.
\end{prop}

\begin{prf}
Let $y \in \mb{R}^n$ be an arbitrary vector.
If $y = 0$, this proposition clearly holds.
If $y \neq 0$, from Definition \ref{def:df}, we obtain
\begin{eqnarray*}
\| y \|_p^*	&=  &	\sup \{ x^T y \: | \: \| x \|_p \leq 1 \} \\
		&\leq&	\sup \{ | x^T y | \: | \: \| x \|_p \leq 1 \} \\
		&\leq&	\sup \biggl\{ \sum_{i=1}^n | x_i | | y_i  | \: | \: \| x \|_p \leq 1 \biggr\} \\
		&\leq&	\max_j | y_j | \left( \sup \biggl\{ \sum_{i=1}^n | x_i | \: | \: \| x \|_p \leq 1 \biggr\} \right) \\
		&= &	\max_j | y_j | \biggl( \sup \{ \| x \|_1 \: | \: \| x \|_p \leq 1 \} \biggr).
\end{eqnarray*}
Since $p \in (0, 1)$, we note that $\|x \|_1 \leq \|x \|_p$ holds \cite{kirmaci2008some}.
Then, we have
\[
\| y \|_p^*	 \leq	\max_j | y_j | \biggl( \sup \{ \| x \|_p \: | \: \| x \|_p \leq 1 \} \biggr)  =  \max_j | y_j | = \|y \|_\infty.
\]
Now, take an arbitrary $i_0 \in \argmax_i | y_i |$,
and define $\bar{x}_i$ as follows:
\[
\bar{x}_i =
\left\{
\begin{array}{ll}
\sign(y_i),	&	\mathrm{if} \quad i = i_0, \\
0,	&	\mathrm{otherwise},
\end{array}
\right.
\]
where
\[
\sign(y_i) =
\left\{
\begin{array}{ll}
1,	&	\mathrm{if} \quad y_i > 0, \\
0,	&	\mathrm{if} \quad y_i = 0, \\
-1,	&	\mathrm{if} \quad y_i < 0.
\end{array}
\right. 
\]
Then, $\|\bar{x} \|_p = 1$ and we have
\[
\|y \|_p^*	= \sup \{ x^T y \: | \: \| x \|_p \leq 1 \}  \geq	\bar{x}^T y = \max_i |y_i | = \|y \|_\infty,
\]
which completes the proof.
\qed
\end{prf}


\bibliographystyle{plain}
\bibliography{reference}

\begin{thebibliography}{10}

\bibitem{alizadeh2003second}
F.~Alizadeh and D.~Goldfarb.
\newblock Second-order cone programming.
\newblock {\em Mathematical programming}, 95(1):3--51, 2003.

\bibitem{ABD2017}
A.~Y. Aravkin, J.~V. Burke, D.~Drusvyatskiy, M.~P. Friedlander, and K.~MacPhee.
\newblock Foundations of gauge and perspective duality.
\newblock {\em arXiv preprint arXiv:1702.08649}, 2017.

\bibitem{caccetta2011globally}
L.~Caccetta, B.~Qu, and G.~Zhou.
\newblock A globally and quadratically convergent method for absolute value
  equations.
\newblock {\em Computational Optimization and Applications}, 48(1):45--58,
  2011.

\bibitem{C2007}
R.~Chartrand.
\newblock Exact reconstruction of sparse signals via nonconvex minimization.
\newblock {\em IEEE Signal Processing Letters}, 14(10):707--710, 2007.

\bibitem{CY2008}
R.~Chartrand and W.~Yin.
\newblock Iteratively reweighted algorithms for compressive sensing.
\newblock In {\em Acoustics, speech and signal processing, 2008. ICASSP 2008.
  IEEE international conference on}, pages 3869--3872. IEEE, 2008.

\bibitem{eldar2009robust}
Y.~C. Eldar and M.~Mishali.
\newblock Robust recovery of signals from a structured union of subspaces.
\newblock {\em IEEE Transactions on Information Theory}, 55(11):5302--5316,
  2009.

\bibitem{F87}
R.~M. Freund.
\newblock Dual gauge programs, with applications to quadratic programming and
  the minimum-norm problem.
\newblock {\em Mathematical Programming}, 38(1):47--67, 1987.

\bibitem{FMP2014}
M.~P. Friedlander, I.~Macedo, and T.~K. Pong.
\newblock Gauge optimization and duality.
\newblock {\em SIAM Journal on Optimization}, 24(4):1999--2022, 2014.

\bibitem{hu2010note}
S.~Hu and Z.~Huang.
\newblock A note on absolute value equations.
\newblock {\em Optimization Letters}, 4(3):417--424, 2010.

\bibitem{hu2011generalized}
S.~Hu, Z.~Huang, and Q.~Zhang.
\newblock A generalized {N}ewton method for absolute value equations associated
  with second order cones.
\newblock {\em Journal of Computational and Applied Mathematics},
  235(5):1490--1501, 2011.

\bibitem{kirmaci2008some}
U.~S. Kirmaci, M.~K. Bakula, M.~E. {\"O}zdemir, and J.~E. Pecaric.
\newblock On some inequalities for $p-$norms.
\newblock {\em Journal of Inequalities in Pure \& Applied Mathematics},
  9(1):1--8, 2008.

\bibitem{luo1996mathematical}
Z.-Q. Luo, J.-S. Pang, and D.~Ralph.
\newblock {\em Mathematical Programs with Equilibrium Constraints}.
\newblock Cambridge University Press, 1996.

\bibitem{mangasarian2007absolute2}
O.~L. Mangasarian.
\newblock Absolute value equation solution via concave minimization.
\newblock {\em Optimization Letters}, 1(1):3--8, 2007.

\bibitem{mangasarian2007absolute}
O.~L. Mangasarian.
\newblock Absolute value programming.
\newblock {\em Computational Optimization and Applications}, 36(1):43--53,
  2007.

\bibitem{mangasarian2009generalized}
O.~L. Mangasarian.
\newblock A generalized {N}ewton method for absolute value equations.
\newblock {\em Optimization Letters}, 3(1):101--108, 2009.

\bibitem{mangasarian2006absolute}
O.~L. Mangasarian and R.~R. Meyer.
\newblock Absolute value equations.
\newblock {\em Linear Algebra and Its Applications}, 419(2):359--367, 2006.

\bibitem{meier2008group}
L.~Meier, S.~van~de Geer, and P.~B{\"u}hlmann.
\newblock The group {L}asso for logistic regression.
\newblock {\em Journal of the Royal Statistical Society: Series B (Statistical
  Methodology)}, 70(1):53--71, 2008.

\bibitem{miao2015generalized}
X.~Miao, J.~Yang, and S.~Hu.
\newblock A generalized {N}ewton method for absolute value equations associated
  with circular cones.
\newblock {\em Applied Mathematics and Computation}, 269:155--168, 2015.

\bibitem{MR2010}
N.~Mourad and J.~P. Reilly.
\newblock Minimizing nonconvex functions for sparse vector reconstruction.
\newblock {\em IEEE Transactions on Signal Processing}, 58(7):3485--3496, 2010.

\bibitem{prokopyev2009equivalent}
O.~Prokopyev.
\newblock On equivalent reformulations for absolute value equations.
\newblock {\em Computational Optimization and Applications}, 44(3):363--372,
  2009.

\bibitem{rohn2004theorem}
J.~Rohn.
\newblock A theorem of the alternatives for the equation $ax+ b|x|= b$.
\newblock {\em Linear and Multilinear Algebra}, 52(6):421--426, 2004.

\bibitem{rohn2009algorithm}
J.~Rohn.
\newblock An algorithm for solving the absolute value equation.
\newblock {\em Electronic Journal of Linear Algebra}, 18(5):589--599, 2009.

\bibitem{stojnic2009reconstruction}
M.~Stojnic, F.~Parvaresh, and B.~Hassibi.
\newblock On the reconstruction of block-sparse signals with an optimal number
  of measurements.
\newblock {\em IEEE Transactions on Signal Processing}, 57(8):3075--3085, 2009.

\bibitem{w2011}
G.~W. Wolf.
\newblock {\em Facility Location: Concepts, Models, Algorithms and Case
  Studies}.
\newblock Taylor \& Francis, 2011.

\bibitem{xue2000efficient}
G.~Xue and Y.~Ye.
\newblock An efficient algorithm for minimizing a sum of $p-$norms.
\newblock {\em SIAM Journal on Optimization}, 10(2):551--579, 2000.

\bibitem{yamanaka2014branch}
S.~Yamanaka and M.~Fukushima.
\newblock A branch-and-bound method for absolute value programs.
\newblock {\em Optimization}, 63(2):305--319, 2014.

\bibitem{yuan2006model}
M.~Yuan and Y.~Lin.
\newblock Model selection and estimation in regression with grouped variables.
\newblock {\em Journal of the Royal Statistical Society: Series B (Statistical
  Methodology)}, 68(1):49--67, 2006.

\bibitem{zhang2009global}
C.~Zhang and Q.~J. Wei.
\newblock Global and finite convergence of a generalized {N}ewton method for
  absolute value equations.
\newblock {\em Journal of Optimization Theory and Applications},
  143(2):391--403, 2009.

\end{thebibliography}

\end{document}